\documentclass[11pt]{article}
\usepackage{bbm}
\usepackage{epsfig}
\usepackage{amssymb,amsmath,amscd, mdwmath}
\usepackage{graphicx,cite,url}

\newcommand{\EOP} { \hfill $\Box$ }

\newtheorem{theorem}{Theorem}[section]
\newtheorem{prop}[theorem]{Proposition}

\begin{document}
\title{\bf Computation of Maxwell's equations on manifold using implicit DEC scheme}

\author{
Zheng Xie$^1$\thanks{E-mail: lenozhengxie@yahoo.com.cn }~~~~Yujie
Ma$^2$\thanks{ E-mail: yjma@mmrc.iss.ac.cn
\  This work is partially
supported by  CPSFFP (No. 20090460102), NKBRPC (No. 2004CB318000),
and NNSFC (No. 10871170) }
\\{\small
$1.$ Center of Mathematical Sciences, Zhejiang University
(310027),China}
\\ {\small $2.$ Key Laboratory of Mathematics Mechanization,}
\\ {\small  Chinese Academy of Sciences,  (100090), China}}

\date{}

\maketitle

\begin{abstract}

  Maxwell's equations can be  solved numerically  in space manifold and the time
   by discrete exterior
calculus as a  kind of lattice gauge theory. Since the stable
conditions of this method is
 very severe restriction, we combine
 the implicit scheme of time variable and
discrete exterior calculus   to derive an unconditional stable
scheme.
 It is an
generation of  implicit Yee-like scheme, since it can be implemented
in space manifold directly.
 The analysis of
its unconditional stability and error is also accomplished.

\end{abstract}


\vskip 0.2cm \noindent {\bf PACS(2010):} 41.20.Jb, 02.30.Jr,
02.40.Sf, 02.60.Cb.

\section{Introduction}

 The  Yee scheme is a commonly employed efficient approach to solve Maxwell's
  equations numerically and so to model wave
propagation problems in the time domain \cite{yee}. Although it is
not a high order method, it is still preferred for many applications
because it preserves important structural features of Maxwell's
equations\cite{bondeson,clemens,nicolet,gross}. Bossavit et al
present the
 Yee-like scheme to
 extend  the Yee scheme to unstructured grids. This  scheme  combines the
 best attributes of the finite element method (unstructured grids)
  and  the Yee scheme (preserving geometric structure)\cite{bossavit1,bossavit2}.
 Stern et al \cite{stern} generalize the Yee-like scheme to
unstructured grids not just in space, but  in 4-dimensional
spacetime by discrete exterior
calculus(DEC)\cite{whitney,arnold,hairer,novikov,auchmann,dimakis,desbrun,meyer,leok,hyman,hiptmair,wise}.
This relaxes the need to take uniform time steps. Based on these
results, we generalize the Yee-like  scheme to
 space manifold and the time, which is  a kind of lattice gauge
 theory\cite{xie-ye-ma}.

   The stable conditions for Yee-like scheme and its generation are
 very severe restriction, and imply that very many time steps will
be necessary to follow the solution over a reasonably large time
interval. In computer simulations of physical processes, explicit
and implicit methods are often used. For some problems,  it takes
much less computational time to use the implicit method with larger
time steps, even taking into account that one needs to solve
equations at each step. Based on these considerations,  Ker\"aen et
al present the unstable conditional implicit
 Yee-like scheme \cite{Janne}. In this paper, we show that    the
implicit scheme  and DEC can be united to find an unconditional
stable scheme (IDEC) for solving Maxwell's equations on space
manifold and the time. The analysis of IDEC's unconditional
stability and error is also accomplished. This scheme  reduces to
the implicit Yee scheme, if choosing rectangular mesh for flat
space, and reduces to the scheme presented  by Ker\"aen et al, if
choosing tetrahedral mesh for flat space.

\section{Preliminaries}

 A discrete differential $k$-form, $k \in \mathbb{Z}$, is the
evaluation   of the differential $k$-form on all $k$-simplices. Dual
forms, i.e., forms  evaluated on the dual cell.  Suppose each
simplex contains its circumcenter. The circumcentric dual cell
$D(\sigma_0)$ of simplex $\sigma_0$ is
 $$ D(\sigma_0):=\bigcup_{\sigma_0\in \sigma_1\in\cdots \in\sigma_r}
 \mathrm{Int}(c(\sigma_0)c(\sigma_1)\cdots c(\sigma_r) ),$$
where   $\sigma_i$ is all the simplices which contains
$\sigma_0$,..., $\sigma_{i-1}$, and $c(\sigma_i)$ is the
circumcenter of $\sigma_i$. In DEC, the   exterior derivative $d$ is
approximated as the transpose of the incidence matrix of $k$-cells
on $k+1$-cells, and the approximated Hodge Star $\ast$ scales the
cells by
  the volumes of the corresponding dual and
primal cells.

 The $2$D or $3$D   space manifold can be approximated by triangles  or tetrahedrons,
  and the time by line segments.    Discrete connection $1-$form or gauge field $A$ assigns to each
element
 in the  set of edges $E$ an element of the gauge group $\mathbb{R}$:
$$A: E\rightarrow \mathbb{R}.$$
Discrete curvature $2-$form is the discrete exterior derivative of
the   discrete  connection $1-$form
$$F=dA: P\rightarrow \mathbb{R}.$$ The value of $F$ on
each element in the set of triangular $P$ is   the coefficient of
Holonomy group of this face. The $2-$form  $F$ automatically
satisfies the discrete Bianchi identity
$$d F=0.\eqno{(1)}$$
For source case, we need discrete current $1-$form  $J$. Let
$A=\sum\limits_EA_i$  and the Lagrangian functional be
 $$\begin{array}{lll}
L(A,J)&=& -\frac{1}{2}\langle dA,  dA\rangle +\langle A,J\rangle,
\end{array}
$$
where
$$\begin{array}{lll}\langle dA,  dA\rangle &:=& (A)_{1\times |E|}(d)_{|E|\times|F|}(*)_{|F|\times|F|}(d)^{T}_{|F|\times|E|}(A)^{T}_{|E|\times 1}\\
\langle A,J\rangle&:=&(A)_{1\times |E|} (*)_{|E|\times|E|}
(J^{T})_{|E|\times 1}.
\end{array}
$$The Hamilton's principle of
stationary action states that this variation must equal zero for any
vary of $A_i$, implying
$$
d^{T}\ast dA=\ast J \eqno{(2)}
$$Since $(d^{T})^2=0$, the discrete
continuity equation can express as:
$$d^{T}\ast J=0.\eqno{(3)}$$ The Eqs.(1-3)  are called
discrete Maxwell's equations, which are invariant under gauge
transformations $A\rightarrow  A+df$ for any $0$-forms.

\section{ IDEC for Maxwell's Equations}
If  allowing for the possibility of magnetic charges and current
discrete $3-$form $\bar{J}$, the symmetric discrete Maxwell's
equations can be written as
$$\begin{array}{lll}dF=\bar{J}&& d^T\ast
F=\ast J, \end{array}\eqno{(4)}$$ with discrete continuity equations
or integrability conditions
$$\begin{array}{lll}
d\bar{J} =0&& d^T\ast J=0. \end{array} \eqno{(5)}$$

\subsection*{Implicit scheme for TE wave}
  The  discrete current  forms,   discrete
curvature $2-$form, and their dual can be written as
$$\begin{array}{lll}
\bar{J}^n =(\rho^n_m,-{J}^{n+\frac{1}{2}}_m\wedge d {t})&& \ast
J^{n}=(-\ast (\rho_e dt)^{n} , \ast
 {J}^{n+\frac{1}{2}}_e)\\
F^{n }=E^{n+1}\wedge d {t}+B^{n} &&\ast F^{n }=H^{n+1 } \wedge d
{t}-D^{n},\end{array}$$
 where $n$ and $n+\frac{1}{2}$   denote the coordinate of the time,
   $E=\sum\limits_{E}{E}_ie^i$ (electric field) is the discrete $1-$form on space,
   $B=\sum
\limits_{P} {B}_i P^i$ (magnetic field) is the discrete $2-$form on
space,  $H=\sum\limits_{P} {H}_i
*P^i$ (magnetizing field) is the dual of $B$   on space,
$D=\sum\limits_{E}  {D}_i
*e^i$(electric displacement field) is the  dual of $E$  on
space,  $\rho_e dt$ (charge density) is the  discrete $1-$form on
time,
   $J_e=\sum\limits_{E}J_{ei}e^i$ (electric current density) is the  discrete $1-$form on space.
   $ \rho_m=\sum\limits_{Tet}\rho_{mi}T^i$(magnetic charges) is the discrete 3-form  on  space,
      ${J}_m=\sum\limits_{P}J_{mi}P^i$
 (current) is the discrete $2-$form on space.

The symmetric discrete Maxwell's equations (4) and (5) can  be
written as
 $$\begin{array}{lll}
d_sB^{n}&=&   \rho^{n }_m\\
d_s E^{n+1}\wedge d {t} &=&-d_{ {t}}B^n-
{J}^{n+\frac{1}{2}}_m\wedge d {t}  \\
d^{T}_sD^{n}&=&\ast (\rho_e dt)^{n}\\
d^{T}_s H^{n+1}\wedge d{t} &=&d^{T}_{{t}}D^{n }+\ast
 {J}^{n+\frac{1}{2}}_e , \end{array}\eqno{(6)}$$
where $d_s$, $d^T_{{s}}$ are the restriction of $d$, $d^T$ on space,
and
$$\begin{array}{lll}d_{ {t}} B^n : = \dfrac{B^{n+1}-B^{n}}{\Delta t} \wedge{d
{t}}&&
   d^{T}_{
{t}}D^{n} := \dfrac{D^{n+1}-D^{n }} {\Delta t}\wedge d
{t}.\end{array} $$

\begin{prop}\label{TH:fundmat} If the initial condition satisfies the first and third equations in
Eqs.(6),  the solution of the second and fourth equations in Eqs.(6)
automatically satisfy Eqs.(6).
\end{prop}
\noindent Proof.  Because the dimension of spacetime is $3+1$ or
$2+1$, therefore
$$d^{T}_s\ast(\rho_e dt)^{n}=0 ~~~~d^{T}_t\ast J^{n+\frac{1}{2}}_e=0~~~~d_s
\rho^n_m=0~~~~d_t  J^{n+\frac{1}{2}}_m\wedge dt=0,$$ and the
continuity equations can be reduced to $$d^T_t\ast( \rho_e
dt)^{n}-d^T_s\ast
J^{n+\frac{1}{2}}_e=0~~~~~-d_sJ^{n+\frac{1}{2}}_m\wedge
dt+d_t\rho^n_m=0.$$ So we have
$$\begin{array}{lll}d^{T}_td^{T}_sD^{n}-d^{T}_t\ast
(\rho_e dt)^{n}&=&-d^{T}_t\ast (\rho_e dt)^{n}-d^{T}_s(d^{T}_s
H^{n+1}\wedge
d{t}-\ast J^{n+\frac{1}{2}}_e)\\
&=&0
\\d_td_sB^{n}-d_t\rho^{n }_m &=& -d_t\rho^{n }_m+d_s(d_s E^{n+1}\wedge d
{t}+ {J}^{n+\frac{1}{2}}_m\wedge d {t}) \\
&=&0.
\end{array}
$$
\EOP

Now we show the implicit scheme (6)  on  the 2D discrete space
manifold and the time.  Take Fig.1  as an example for a part of 2D
space mesh, in which $e_1$,..., $e_5$ are edges,    $P_1$, $P_2$ are
 triangles,    $\ast e_1$ is the dual of
$e_1$. The second and fourth equations in Eqs.(6) based on Fig.1 are
\begin{equation*}
\left.
    \begin{aligned}
&\frac{ {D}^{n+1}_1 -  {D}^{n }_1 }{\Delta
t}+{J^{n+\frac{1}{2}}_{e1}}
=\frac{ {H}^{n+1}_{1} - {H}^{n+1}_{2} }{|*e_1|}&\\
&-\frac{ {B}^{n+1}_1 - {B}^{n}_1 } {\Delta
t}-{J}^{n+\frac{1}{2}}_{m1}=\frac { {E}^{n+1}_1 |e_1|+ {E}^{n+1}_2
|e_2|+ {E}^{n+1}_3 |e_3|}{|P_1|}.&
\end{aligned}
  \right\} \eqno{(7)}
\end{equation*}
 where $|~|$ denotes the measure of forms and dual.
The summation on the right is orient, that is to say, inverse the
orientation of $e_i$, then multiply $-1$ with ${E}_i$. Notice that a
significant difference from Yee-like scheme is that $\ast e_i$ is a
polyline. Eqs.(7) can be implemented on 2D discrete space manifold
directly (see Fig.3) so is  a generation of implicit Yee-like
scheme.
$$
\begin{minipage}{0.99\textwidth}
\begin{center}\includegraphics[scale=0.4]{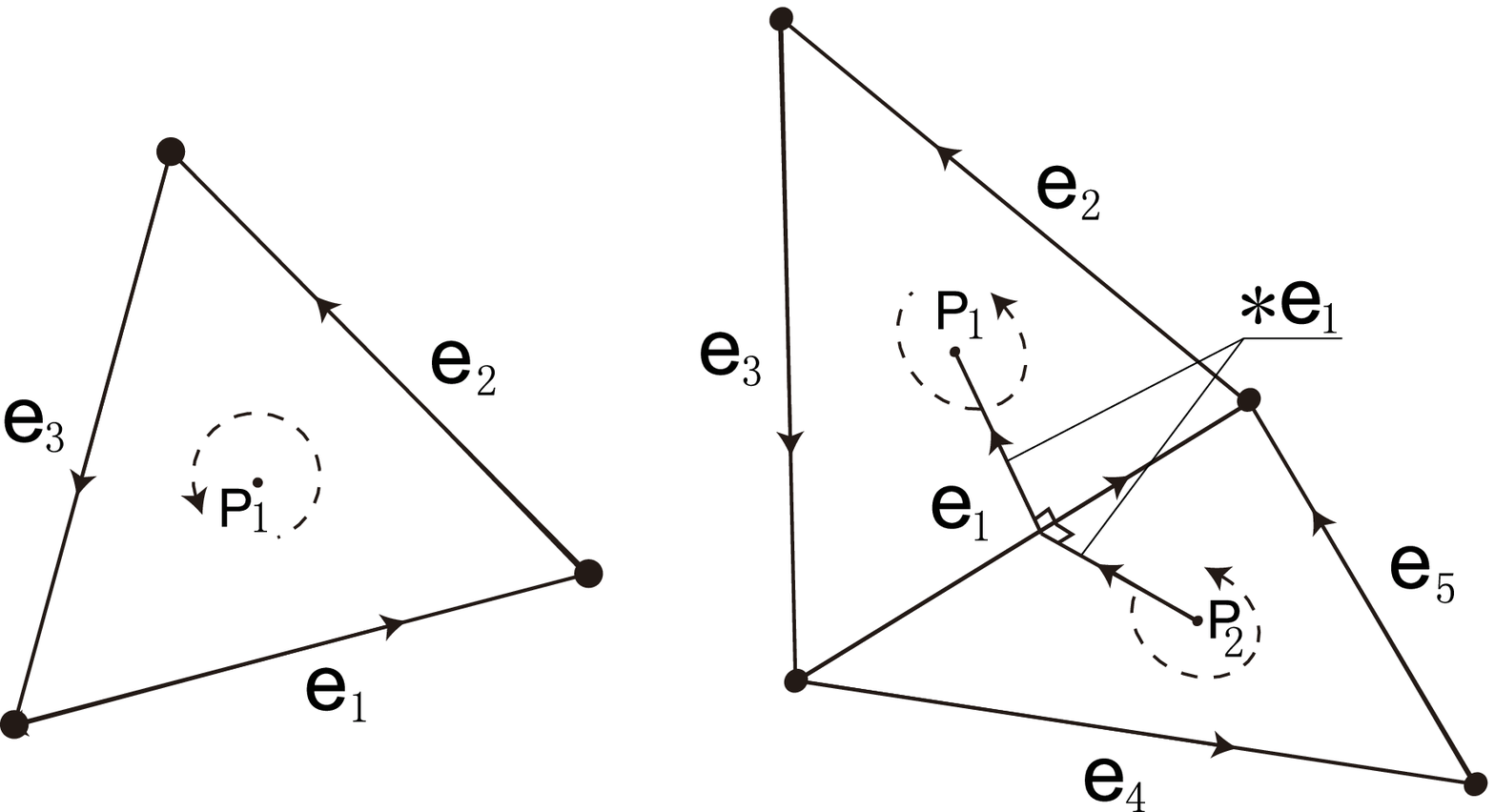}
\end{center}
\begin{center}
{Figure 1: edge and face with  direction  }
\end{center}
\end{minipage}
$$

 In the
absence of magnetic or dielectric materials, there are relations
$$\begin{array}{lll}D_i=\varepsilon_0 E_i &&
B_i=\mu_0H_i,\end{array}\eqno{(8)}$$ where $\varepsilon_0$ and
$\mu_0$ are two universal constants, called the permittivity of free
space and permeability of free space, respectively. With relations
(8), Eqs.(7) can be rewritten
  into an implicit scheme (9) for
  TE  wave.
\begin{equation*}
\left.
    \begin{aligned}
&\varepsilon_0\dfrac{ {E}^{n+1}_1 - {E}^{n }_1 } {\Delta t}
+J^{n+\frac{1}{2}}_{e1}=\dfrac{ {H}^{n+1}_{1} - {H}^{n+1}_{2} }{|*e_1|}&\\
&\mu_0\dfrac{ {H}^{n+1}_1-  {H}^{n}_1 }{\Delta
t}+J^{n+\frac{1}{2}}_{m1} =-\dfrac{ {E}^{n+1}_1 |e_1|+ {E}^{n+1}_2
|e_2|+ {E}^{n+1}_3 |e_3|}{|P_1|}&
\end{aligned}
  \right\}\mathrm{TE}\eqno{(9) }
\end{equation*}

\subsection*{Implicit scheme  for TM wave}

If writing
$$\begin{array}{lll}
F^{n}= H^{n+1}\wedge d {t}-D^{n}  &&\ast
F^{n}=-E^{n+1 } \wedge d {t}-B^{n}\\
 \bar{J}^n=(-\rho^n_e,{J}^{n+\frac{1}{2}}_e\wedge dt)&&
  \ast{J}=(-\ast(\rho_m dt)^{n},\ast{J}^{n+\frac{1}{2}}_m),
\end{array}$$
 where  $H=\sum\limits_{E} {H}_i
e^i$  is the discrete   $1-$form   on space, $D=\sum\limits_{P}
{H}_i P^i$  is the discrete   $2-$form  on space,
  $E=\sum\limits_{P} {D}_i
*P^i$  is the dual of $D$  on space,
    $B=\sum\limits_{E} {B}_i
*e^i$   is the dual of $H$  on space,
 $\rho_e=\sum\limits_{Tet}\rho_{ei}T^i$  is the discrete   $3-$form
on space,  ${J}_e=\sum\limits_{P}J_{ei}P^i$  is the discrete
$2-$form on space,  $\rho_m dt$  is the discrete  $1-$form   on
time,  ${J}_m=\sum\limits_{E}J_{mi}E^i$  is the discrete   $1-$form
on space,  the discrete Maxwell's equations can  be rewritten as
$$\begin{array}{lll}
d_sD^{n}&=& \rho^n_e\\
d_s H^{n+1}\wedge d {t} &=& d_{ {t}}D^n+J^{n+\frac{1}{2}}_e\wedge dt\\
d^{T}_sB^{n }&=&\ast (\rho_m  dt)^{n}\\
d^{T}_s E^{n+1}\wedge d{t} &=&-d^{T}_{{t}}B^{n}-\ast
{J}^{n+\frac{1}{2}}_m  .
\end{array}\eqno{(10)}$$
\begin{prop}\label{TH:fundmat} If the initial condition satisfies the first and third equations in
Eqs.(10),  the solution of the second and fourth equations in
Eqs.(10) automatically satisfy Eqs.(10).
\end{prop}
\noindent Proof. Because the dimension of spacetime is $3+1$ or
$2+1$, therefore $$d^{T}_s\ast(\rho_m dt)^{n}=0~~~d^{T}_t\ast
J^{n+\frac{1}{2}}_m=0~~~d_s \rho^{n}_e=0~~~d_t
J^{n+\frac{1}{2}}_e\wedge dt=0,$$ and the continuity equations can
be reduced to$$d^{T}_t\ast( \rho_m dt)^{n}-d^{T}_s(\ast
J^{n+\frac{1}{2}}_m)=0~~~~~d_sJ^{n+\frac{1}{2}}_e\wedge
dt-d_t\rho^n_e=0.$$ So we have
$$\begin{array}{lll}d_td_sD^{n}-d_t
\rho_e^{n}&=&-d_t \rho_e^{n}-d_s(d_s H^{n+1}\wedge
d{t} - J^{n+\frac{1}{2}}_e\wedge dt) \\
&=&0
\\d^T_td^T_sB^{n }-d^T_t\ast (\rho_m  dt)^{n }  &=&
-d^T_t\ast (\rho_m  dt)^{n}+d^T_s(d^T_s E^{n+1}\wedge d
{t}+ \ast {J}^{n+\frac{1}{2}}_m )\\
&=&0.
\end{array}
$$
\EOP

 Now we show the scheme (10)  on  the  2D discrete space
manifold and the time.  The second and fourth equations in Eqs.(10)
based on Fig.1 are
\begin{equation*}
\left.
    \begin{aligned}
&\frac{ {B}^{n+1}_1 -  {B}^{n }_1 }{\Delta t}+{J^{n+
\frac{1}{2}}_{m1}}
=-\frac{ {E}^{n+1}_{1} - {E}^{n+1}_{2} }{|*e_1|}&\\
& \frac{ {D}^{n+1}_1 - {D}^{n}_1 } {\Delta
t}+{J^{n+\frac{1}{2}}_{e1}}=\frac { {H}^{n+1}_1 |e_1|+ {H}^{n+1}_2
|e_2|+ {H}^{n+1}_3 |e_3|}{|P_1|}.&
\end{aligned}
  \right\}\eqno{(11)}
\end{equation*}
 With relations (8),
Eqs.(11) can be rewritten
  into an implicit scheme (12) for
  TM  wave.
\begin{equation*}
\left.
 \begin{aligned}
    &\varepsilon_0\dfrac{ {E}^{n+1}_1- {E}^n_1 }
{\Delta t}+ {J}^{n+\frac{1}{2}}_{e1} =\dfrac{ {H}^{n+1}_1|e_1|+
{H}^{n+1}_2
|e_2|+ {H}^{n+1}_3|e_3|}{|P_1|}&\\
&\mu_0\dfrac{ {H}^{n+1}_1 -  {H}^{n }_1}{\Delta t}+
{J}^{n+\frac{1}{2}}_{m1} =-\dfrac{ {E}^{n+1}_{1} -
{E}^{n+1}_{2}}{|*e_1|} &\end{aligned}
  \right\}\mathrm{TM}\eqno{(12)}
\end{equation*}

\subsection*{General schemes }

For real world materials, the constitutive relations are not simple
proportionalities, except approximately. The relations can usually
still be written:
$$\begin{array}{lll}D=\varepsilon  E &&
B=\mu H,\end{array}\eqno{}$$ but $\varepsilon$ and $\mu$ are not, in
general, simple constants, but rather functions. With Ohm's law
$$E=\dfrac{1}{\sigma}J ~~~~~~J_m=\dfrac{1}{\sigma_m}H ,$$where $\sigma$ is the electrical conductivity and $\sigma_m$ is
magnetic conductivity. The IDEC schemes can be written as
\begin{equation*}
\left.
    \begin{aligned}
&\varepsilon\dfrac{ {E}^{n+1}_1 - {E}^{n }_1 } {\Delta
t}+{\sigma}\dfrac{ {E}^{n+1}_1 + {E}^{n }_1 }{2}
=\dfrac{ {H}^{n+1}_{1} - {H}^{n+1}_{2} }{|*e_1|}&\\
&\mu\dfrac{ {H}^{n+1}_1-  {H}^{n}_1 }{\Delta t}+\sigma_m\dfrac{
{H}^{n+1}_1 + {H}^{n}_1} {2} =-\dfrac{ {E}^{n+1}_1 |e_1|+
{E}^{n+1}_2 |e_2|+ {E}^{n+1}_3 |e_3|}{|P_1|},&
\end{aligned}
  \right\}\mathrm{TE}
\end{equation*}
\begin{equation*}
\left.
 \begin{aligned}
    &\varepsilon\dfrac{ {E}^{n+1}_1- {E}^{n }_1 }
{\Delta t}+{\sigma}\dfrac{ {E}^{n+1}_1 + {E}^{n }_1 }{2} =\dfrac{
{H}^{n+1}_1
|e_1|+ {H}^{n+1}_2 |e_2|+ {H}^{n+1}_3 |e_3|}{|P_1|}~~~~~~&\\
&\mu\dfrac{ {H}^{n+1}_1 -  {H}^{n }_1 }{\Delta t}+\sigma_m\dfrac{
{H}^{n+1}_1 + {H}^{n }_1 }{2} =-\dfrac{ {E}^{n+1}_{1} -
{E}^{n+1}_{2} }{|*e_1|}.&
\end{aligned}
  \right\}\mathrm{TM}
\end{equation*}

\section{Stability, convergence  and accuracy}

Now, we analyze the stability for scheme (9). The analysis for
scheme (12) can be done in the same way. Suppose the fields to be:
$$\begin{array}{lll}E^{n+1}_i&=&E^n_i\xi\\
H^{n+1}_2&=&H^n_1 \cos(k|\ast e_1|)\xi,
\end{array}\eqno{(13)}
$$
where $\xi$ is the growth factor of time, and $k$ is the spatial
frequency spectrum. Substituting (13) into scheme (9), we obtain
$$\begin{array}{lll}E^n_1\xi&=&
E^n_1+\dfrac{\Delta t}{\varepsilon |\ast e_1|}(1-\cos(k|\ast
e_1|))H^n_1\xi\\
H^n_1\xi&=&H^n_1-\dfrac{\Delta
t}{\mu|P_1|}\left(E^n_1|e_1|+E^n_2|e_2|+E^n_3|e_3|  \right)\xi
\end{array}\eqno{(14)}$$
Rewrite the first equation of Eqs.(14) as
$$E^n_1=\dfrac{\Delta t}{\varepsilon |\ast e_1|(\xi-1)}(1-\cos(k|\ast e_1|))H^n_1\xi,$$
and substitute it into the second equation of Eqs.(14) to obtain
$$\begin{array}{lll}H^n_1\xi&=&H^n_1-\dfrac{(c\Delta t)^2}{ |P_1|}\left((1-\cos(k|\ast e_1|))H^n_1\dfrac{|e_1|}{|\ast e_1|}
+(1-\cos(k|\ast e_2|))H^n_1\dfrac{|e_2|}{|\ast
e_2|}\right.\\
&&\left.+(1-\cos(k|\ast e_3|))H^n_1\dfrac{|e_3|}{|\ast e_3|} \right)
\dfrac{\xi^2}{\xi-1}
\end{array}$$
Therefore, we obtain a quadratic equation for $\xi$ as follows:
$$ (1+M)\xi^2-2\xi+1=0,\eqno{(15)}$$
where $$\begin{array}{lll}M&=&\dfrac{(c\Delta
t)^2}{|P_1|}\left((1-\cos(k|\ast e_1|))\phi_1\dfrac{|e_1|}{|\ast
e_1|}
+(1-\cos(k|\ast e_2|))\phi_1\dfrac{|e_2|}{|\ast e_2|}\right.\\
&&\left. +(1-\cos(k|\ast e_3|))\phi_1\dfrac{|e_3|}{|\ast e_3|}
\right)\geq 0
\end{array}
$$
The discriminant of Eq.(15) is
$$4-4(1+M)=-4M\leq0.$$
So
$$|\xi|=\frac{1}{\sqrt{1+M}}\leq1.$$
That is to say scheme (9) is unconditional stability.

By the definition of truncation error, the exact solution
  of Maxwell's equations satisfy the same relation as IDEC scheme
  except for an additional term $O((\Delta t)^2+\Delta t|\ast e|)$.
   This expresses the consistency, and so  convergence for  IDEC scheme by  Lax equivalence theorem
   (consistency $+$ stability $=$ convergence).   The
derivative of IDEC scheme  are approximated by first order
difference. Equivalently, $H$ and $E$ are approximated by linear
interpolation functions. Consulting the
  definition about accuracy of finite volume method, we
can also say that IDEC scheme    has first order temporal and
spacial accuracy.

\section{Implementation}

The IDEC scheme of Maxwell's equations was implemented  on C++
platform, consisting of the following steps:
\begin{itemize}
  \item [1.]Set the simulation parameters. These
are the dimensions of the computational grid and the size of the
time step, etc.

  \item [2.] Set the propagating media parameters.

 \item [3.] Initialize the mesh indexes.

  \item [4.]Assign current transmitted signal.

  \item [5.]Compute the value of all spatial nodes  and temporarily store the
  result  in the circular buffer for further computation.

  \item [6.]Visualize the currently computed grid of spatial nodes.

   \item [7.]Repeat the whole process  from the step 4, until reach the desired time.
\end{itemize}
 Fig.3 and Fig.4 exhibit the propagation of Gaussian pulse  on rabbit and
 sphere simulated by IDEC.
 $$
\begin{minipage}{0.99\textwidth}
\begin{center}\includegraphics[scale=0.45]{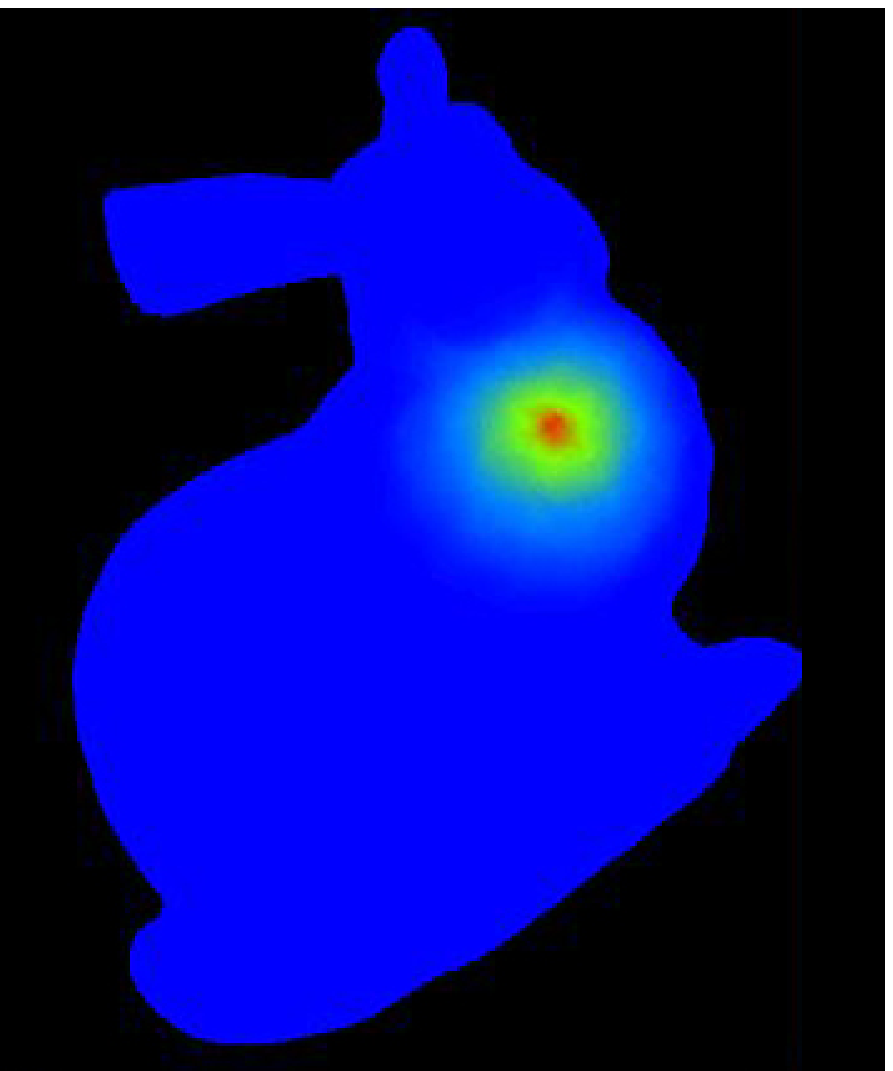}
\includegraphics[scale=0.45]{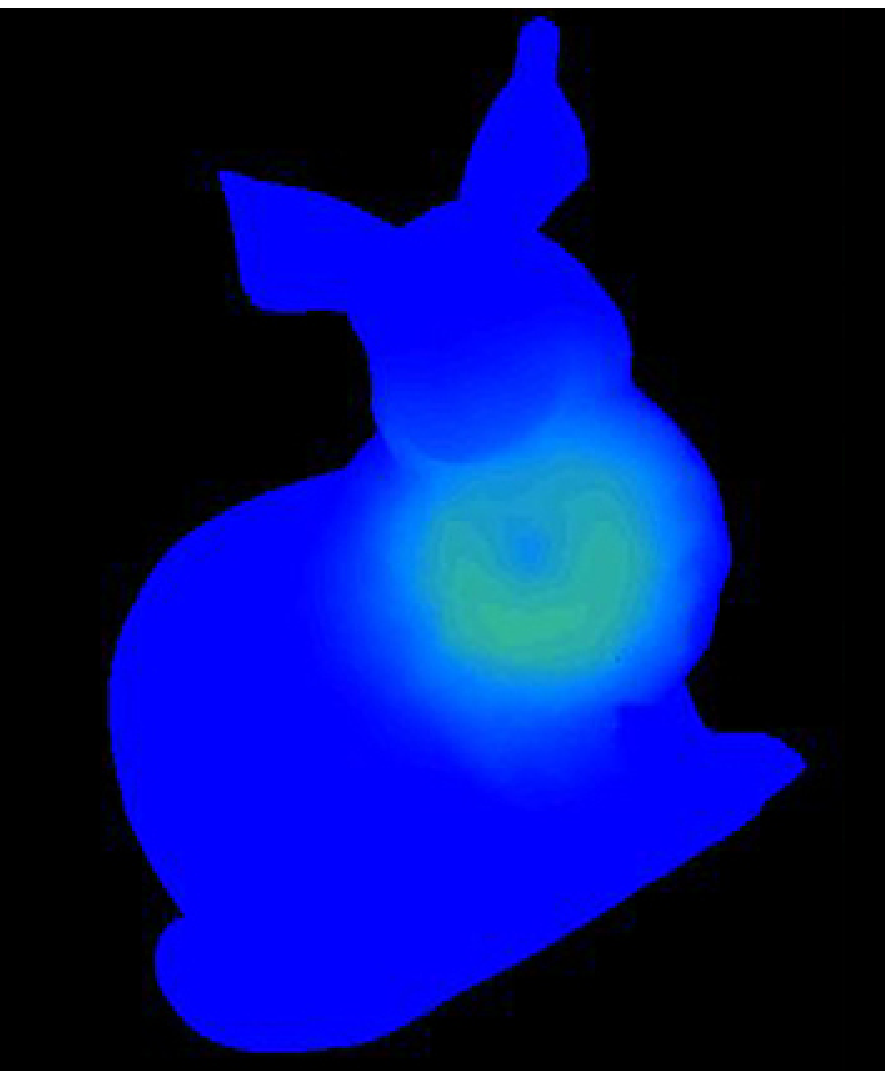}
\includegraphics[scale=0.45]{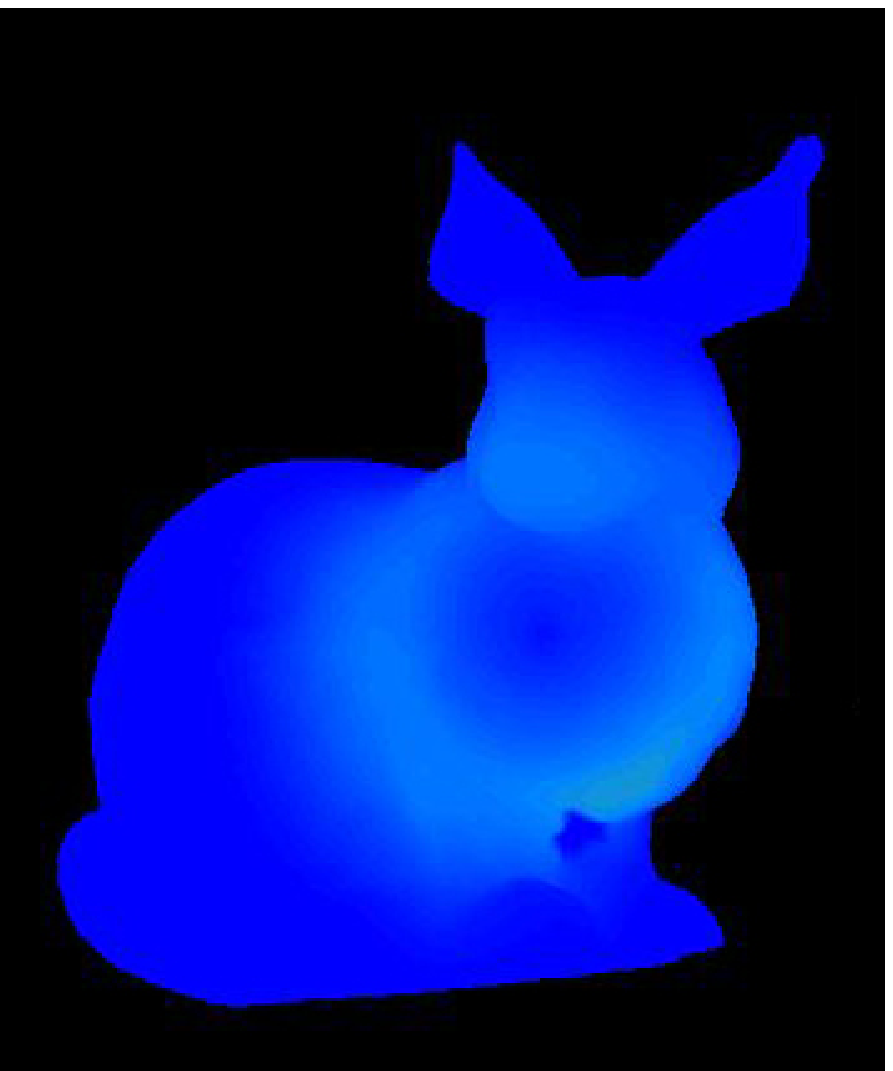}
\end{center}
\begin{center}
{Figure  4 }
\end{center}
\end{minipage}
$$
$$
\begin{minipage}{0.99\textwidth}
\begin{center}\includegraphics[scale=0.5]{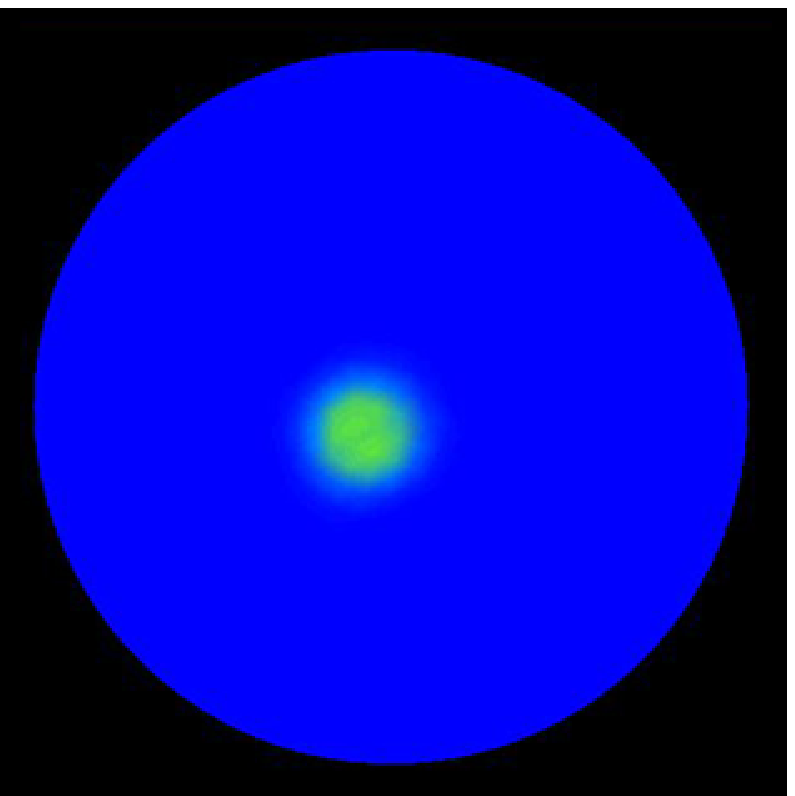}
\includegraphics[scale=0.5]{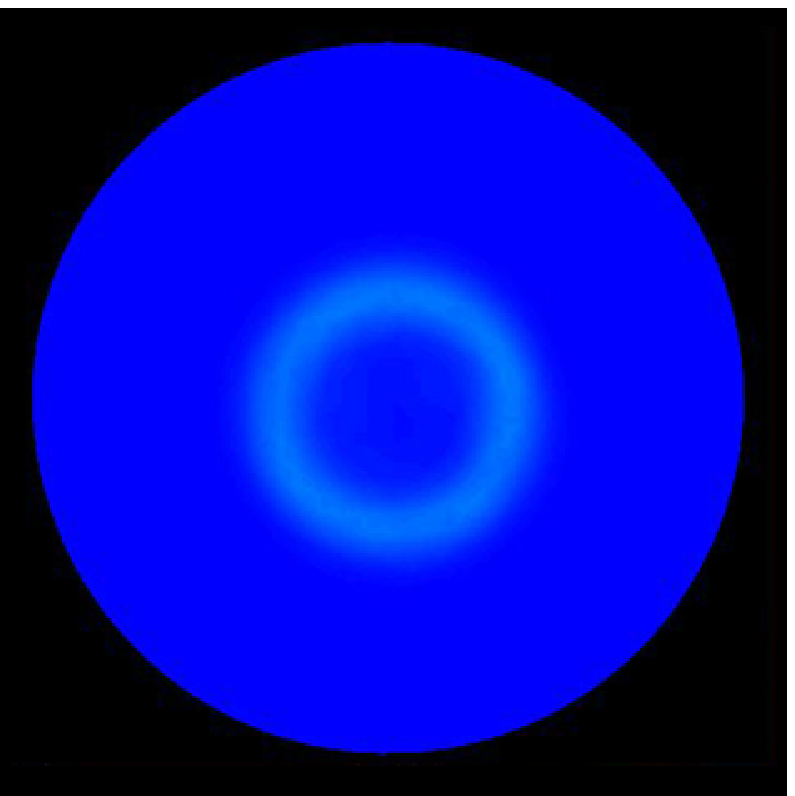}
\includegraphics[scale=0.5]{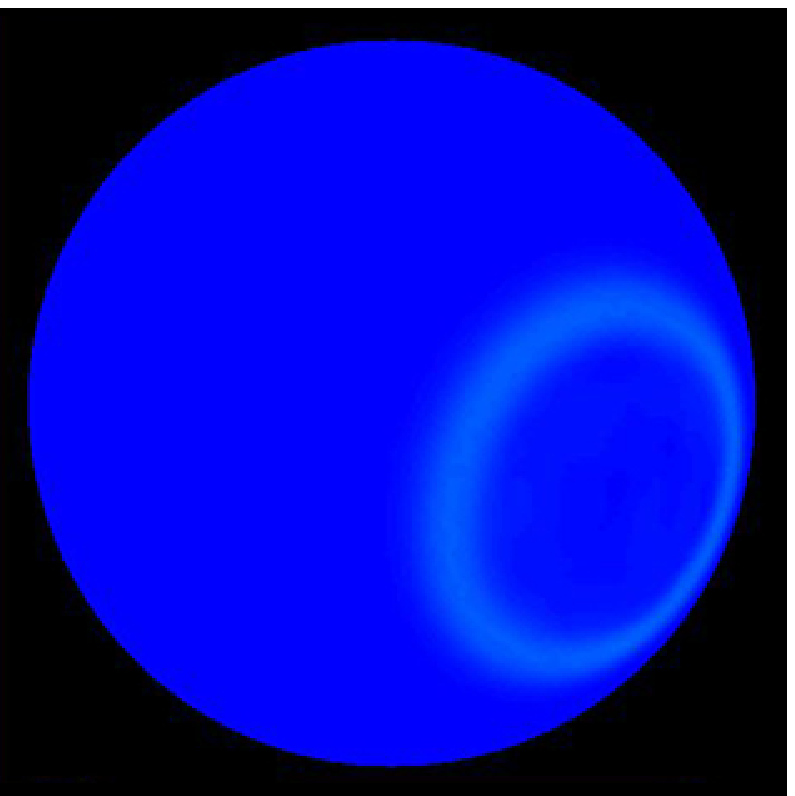}
\end{center}
\begin{center}
{Figure 5 }
\end{center}
\end{minipage}
$$

\newpage

\end{document}